\documentclass[a4paper,12pt,reqno]{amsart}%
\usepackage{eurosym}
\usepackage{amssymb}
\usepackage{amsfonts}
\usepackage{amsmath}
\usepackage{amsthm}
\usepackage{graphicx}
\usepackage[numbers,sort&compress]{natbib}
\usepackage[bookmarks=true,colorlinks=true,linkcolor=blue,citecolor=blue]%
{hyperref}%
\providecommand{\U}[1]{\protect\rule{.1in}{.1in}}
\newtheorem{theorem}{Theorem}
\newtheorem{remark}{Remark}[section]

\begin{document}

\begin{center}
{\Large \textbf{Weighted Gaussian Approximations for Increments of the Uniform
Empirical and Quantile Processes: Fixed-Endpoint Extensions to the
Finite-Count Scale}}\medskip\medskip

{\large Abdelhakim Necir}\medskip

{\small \textit{Laboratory of Applied Mathematics, Mohamed Khider University,
Biskra, Algeria}}\medskip\medskip

{\small Email address: \texttt{ah.necir@univ-biskra.dz}}
\end{center}

\smallskip\smallskip\smallskip\noindent\textbf{Abstract.} We establish
weighted Gaussian approximations for the uniform empirical and quantile
processes and for their increments ending at a fixed point $t\in(0,1)$. We
first place the classical weighted approximations for the ordinary processes
in a common framework and then show that the corresponding increment
approximations remain valid uniformly down to the finite-count scale
$\lambda/n$, for every fixed $\lambda>0$. For the empirical increments, the
proof splits the sample at $t$, couples the two resulting conditional
empirical processes with independent Brownian bridges, and approximates the
binomial fluctuation at $t$ by a Gaussian variable. The three Gaussian
components are then combined into a single standard Brownian bridge. For the
quantile increments, the R\'enyi representation and a reversal of the relevant
exponential spacings reduce the problem to the weighted approximation of an
ordinary uniform quantile process. The resulting bounds hold for $0\leq
\nu<1/4$ in the empirical case and for $0\leq\eta<1/2$ in the quantile case.
As an application of the empirical increment approximation, we derive
simultaneous weighted Gaussian approximations for the censored and uncensored
empirical subdistribution-tail processes arising under random right
censoring.\medskip

\medskip

\newpage

\section{Introduction\label{Sec1}}

\noindent Gaussian approximations for empirical and quantile processes
constitute a fundamental topic in probability theory and mathematical
statistics; see, for instance, \cite{SW1982,CsCsHM86}. Beyond their intrinsic
probabilistic interest, they provide powerful tools for studying statistics
that can be represented as functionals of empirical, quantile, tail empirical,
or related stochastic processes.

\smallskip

\noindent Such approximations play a particularly important role in
extreme-value theory, where many estimators and test statistics are
constructed from intermediate or upper order statistics. They were used by
\cite{cdm85} in the study of kernel estimators of the tail index, by
\cite{CSM1985,CSM1986,CSHM1986} in establishing limit theorems for sums of
extreme values, and by \cite{peng2001} in the estimation of the mean of
heavy-tailed distributions. These contributions illustrate the central role of
Gaussian process approximations in the asymptotic analysis of extreme-value statistics.

\smallskip

\noindent Under random right censoring, the probabilistic structure becomes
more involved. Each observation consists of the minimum of the variable of
interest and an independent censoring variable, together with an indicator
specifying whether the variable of interest has been observed. Tail-related
statistics are therefore naturally expressed in terms of empirical
subdistribution processes associated with uncensored and censored
observations. Gaussian approximations have been used in this framework by
\cite{BMN2015} to study censored extreme-value index estimators, by
\cite{MNS2025} to investigate tail product-limit processes, and by
\cite{GNM2026A,GNM2026B} to develop kernel-based procedures for heavy-tailed
censored data.

\smallskip

\noindent There is a somewhat unusual history behind the increment results
established in this paper. In some of our recent applications, the Gaussian
approximation proved here was used as a natural consequence inferred from the
existing strong approximation theory, although a complete and independent
proof at the finite-count scale was not available at that time. This inference
was not merely heuristic. It was rooted in the finite-dimensional theory of
increments: for a fixed increment length, the law of a uniform empirical
increment depends only on the length of the underlying interval and not on its
starting point. Moreover, since $t$ is fixed in the interior of $(0,1)$,
factors involving the distance from the right endpoint play no essential
asymptotic role. These observations naturally suggested the form of the
approximation. The nontrivial step is to turn this insight into a weighted
uniform approximation over the full range of increment lengths, down to
$\lambda/n$. The strong approximation results of \cite{CsCsHM86} provide the
essential probabilistic tools, whereas the present work makes explicit the
transfer from those tools to the empirical and quantile increment processes.

\noindent Random censoring also provides a natural setting in which increments
of empirical processes arise. Under the uniform representation of
\cite{EnKo92}, one empirical subdistribution component is represented by an
ordinary uniform empirical process, whereas the other is represented by an
increment of the same process. Since both components enter product-limit
estimators and related extreme-value procedures, their joint asymptotic
analysis requires Gaussian approximations that preserve this common dependence structure.

\smallskip

\noindent We now introduce the processes considered throughout the paper. Let
\[
U_{n}(u) := \frac{1}{n}\sum_{i=1}^{n}\mathbb{I}\{U_{i}\leq u\}, \qquad0\leq
u\leq1,
\]
be the empirical distribution function based on i.i.d.\ random variables
$U_{1},U_{2},\ldots$, uniformly distributed on $(0,1)$, where $\mathbb{I}%
\{A\}$ denotes the indicator of the event $A$. The associated uniform
empirical process is
\[
\alpha_{n}(u) := \sqrt n\{U_{n}(u)-u\}, \qquad0\leq u\leq1.
\]
Let
\[
V_{n}(u) := \inf\left\{  v\in[0,1]:U_{n}(v)\geq u\right\}  , \qquad0<u\leq1,
\]
with $V_{n}(0):=0$, denote the generalized inverse of $U_{n}$. Equivalently,
\[
V_{n}(u)=U_{k:n} \quad\text{whenever}\quad\frac{k-1}{n}<u\leq\frac{k}{n},
\qquad k=1,\ldots,n,
\]
where $U_{1:n}\leq\cdots\leq U_{n:n}$ are the corresponding uniform order
statistics. The uniform quantile process is defined by
\[
\beta_{n}(u) := \sqrt n\{u-V_{n}(u)\}, \qquad0\leq u\leq1.
\]

\smallskip

\noindent For any function $f$ defined on $[0,1]$, set
\[
f(s;t) := f(t)-f(t-s), \qquad0\leq s\leq t\leq1.
\]
In particular,
\begin{align*}
\alpha_{n}(s;t)  &  = \alpha_{n}(t)-\alpha_{n}(t-s)\\
&  = \sqrt n\{U_{n}(t)-U_{n}(t-s)-s\}\\
&  = \frac{1}{\sqrt n}\sum_{i=1}^{n} \left\{  \mathbb{I}\{t-s<U_{i}\leq
t\}-s\right\}  , \qquad0\leq s\leq t,
\end{align*}
and
\[
\beta_{n}(s;t) = \beta_{n}(t)-\beta_{n}(t-s) = \sqrt n\left\{  s-\left(
V_{n}(t)-V_{n}(t-s)\right)  \right\}  , \qquad0\leq s\leq t.
\]
Thus, $\alpha_{n}(s;t)$ is the centered and normalized number of sample
observations falling in the interval $(t-s,t]$, whereas $\beta_{n}(s;t)$ is
the centered and normalized increment of the uniform empirical quantile
function. On the empirical-quantile grid, the latter is expressed in terms of
spacings between uniform order statistics, a fact that will be exploited
through the R\'enyi representation.

\smallskip

\noindent Gaussian approximations for such increments have been studied in
several forms. As a consequence of Theorem~1.2 of \cite{SW1982}, there exist
versions of the uniform empirical processes and a sequence of standard
Brownian bridges $\{B_{n}^{\mathrm{SW}}\}_{n\geq1}$, defined on a common
probability space, such that, for every fixed $t\in(0,1)$, every $c>0$, and
every $0\leq\eta<1/2$,
\begin{equation}
\sup_{c(\log n)/n\leq s<t} \frac{\left\vert \alpha_{n}(s;t)-B_{n}%
^{\mathrm{SW}}(s;t)\right\vert } {s^{\eta}} =o_{\mathbb{P}}(1).
\label{SW-emp-incr}%
\end{equation}
\cite[Theorem~4.6.1]{CsCsHM86} considered the corresponding global increment
problem for the uniform quantile process, uniformly over all intervals whose
lengths exceed a logarithmic threshold. More precisely, writing
\[
h:=b-a,
\]
they treated simultaneously all $a$ and $h$ satisfying
\[
c\frac{\log n}{n}\leq h\leq1,\qquad0\leq a\leq1-h.
\]
For an appropriate sequence of standard Brownian bridges $\{B_{n}^{\mathrm{C}%
}\}_{n\geq1}$, put
\[
D_{n}^{\left(  q\right)  }(a,h) := \beta_{n}(a+h)-\beta_{n}(a) -\left\{
B_{n}^{\mathrm{C}}(a+h)-B_{n}^{\mathrm{C}}(a)\right\}  .
\]
Then, for every $0\leq\eta<1/2$,
\[
\sup_{c(\log n)/n\leq h\leq1} \sup_{0\leq a\leq1-h} \frac{\left\vert
D_{n}^{\left(  q\right)  }(a,h)\right\vert }{h^{\eta}} =o(1), \qquad
\text{a.s.}
\]
In Remark~4.6.1, they observed that the analogous empirical-process result
follows from an appropriate Koml\'{o}s--Major--Tusn\'{a}dy construction.
Writing
\[
D_{n}^{\left(  e\right)  }(a,h) := \alpha_{n}(a+h)-\alpha_{n}(a) -\left\{
B_{n}^{\mathrm{C}}(a+h)-B_{n}^{\mathrm{C}}(a)\right\}  ,
\]
the corresponding assertion takes the form
\[
\sup_{c(\log n)/n\leq h\leq1} \sup_{0\leq a\leq1-h} \frac{\left\vert
D_{n}^{\left(  e\right)  }(a,h)\right\vert }{h^{\eta}} =o(1), \qquad
\text{a.s.}
\]
Thus, the lower bound $c(\log n)/n$ is sufficient when all initial points and
all admissible increment lengths are considered simultaneously.

\smallskip

\noindent For increments ending at a fixed terminal point $t$, the preceding
formulation corresponds to taking
\[
a=t-s\qquad\text{and}\qquad h=s.
\]
It then yields an approximation over the range
\[
s\geq c\frac{\log n}{n}.
\]
This logarithmic restriction, however, need not be imposed for the empirical
increment process when the terminal point $t$ is fixed. The weighted central
limit theorem of \cite{Alex87}, specialized to intervals ending at $t$ and
combined with a suitable representation argument, implies that versions of the
empirical processes and a sequence of standard Brownian bridges $\{B_{n}%
^{\mathrm{A}}\}_{n\geq1}$ may be constructed on a common probability space
such that, for every fixed $t\in(0,1)$, every fixed $\lambda>0$, and every
$0\leq\eta<1/2$,
\begin{equation}
\sup_{\lambda/n\leq s<t} \frac{\left\vert \alpha_{n}(s;t)-B_{n}^{\mathrm{A}%
}(s;t)\right\vert } {s^{\eta}} =o_{\mathbb{P}}(1). \label{Alex-emp-incr}%
\end{equation}
Thus, for empirical increments ending at a fixed point, the finer scale
$\lambda/n$ is already attainable.

\smallskip

\noindent This scale corresponds to a finite-count regime. Indeed, an interval
of length $s$ contains, on average, $ns$ observations. When $s=\lambda/n$, the
expected number of observations contributing to the increment is exactly
$\lambda$. At this level, the discrete nature of the empirical process and the
possible contribution of individual sample observations can no longer be neglected.

\smallskip

\noindent The fixed-terminal-point setting also reveals the structural reason
why the scale $\lambda/n$ is natural. For each fixed $t$, the empirical
increment $\alpha_{n}(s;t)$ is, after a rotation of the uniform sample, an
ordinary uniform empirical process indexed by the interval length $s$.
Similarly, on the empirical-quantile grid, the quantile increment $\beta
_{n}(s;t)$ can be expressed, through the R\'enyi representation and a reversal
of the relevant exponential spacings, as an ordinary uniform quantile process.
In both cases, the starting point $t$ disappears from the finite-dimensional
law, while the genuine difficulty lies in obtaining a weighted approximation
uniformly over all $\lambda/n\leq s<t$.

\smallskip

\noindent The contribution of the present paper is twofold. First, for each
fixed terminal point $t\in(0,1)$, we construct a single sequence of standard
Brownian bridges that simultaneously approximates the ordinary uniform
empirical process and its increment process. The resulting weighted
approximation bounds hold uniformly down to the scale $\lambda/n$, for every
fixed $\lambda>0$ and every $0\leq\nu<1/4$. Second, we establish the
corresponding weighted approximation for the uniform quantile increment
process at the same finite-count scale. The latter result is obtained by
combining the R\'enyi representation with the reversal of the exponential
spacings associated with the fixed terminal point. It therefore extends, in
the fixed-$t$ setting, the logarithmic range appearing in \cite[Theorem~4.6.1]%
{CsCsHM86} to the finer range beginning at $\lambda/n$. The precise statements
are given in Theorem~\ref{Theo1}.

\smallskip

\noindent The ordinary empirical- and quantile-process approximations
appearing in Theorem~\ref{Theo1}, when considered separately, belong to the
classical weighted approximation theory. For the empirical part, the essential
point is that the ordinary process and its fixed-terminal-point increment
process are realized with the same sequence of Brownian bridges. For the
quantile part, the reversal argument identifies the increment process with an
ordinary uniform quantile process and thereby makes the classical weighted
approximation available down to the finite-count scale.

\smallskip

\noindent The use of the same Brownian bridge in the two empirical
approximations is essential for the censoring application. It preserves the
joint Gaussian dependence between the ordinary empirical component and the
increment component appearing in the uniform representation of the empirical
subdistributions. Moreover, the coupling may be constructed independently of
the particular values of $\lambda$ and $\nu$, although it may depend on the
fixed terminal point $t$.

\smallskip

\noindent A direct subtraction of two ordinary-process approximations does not
provide the weighted increment bound stated in Theorem~\ref{Theo1}. Indeed, an
ordinary Koml\'{o}s--Major--Tusn\'{a}dy remainder of order $n^{-1/2}\log n$
would yield, at the smallest admissible scale $s=\lambda/n$,
\[
\frac{n^{\nu}n^{-1/2}\log n}{(\lambda/n)^{1/2-\nu}} = \lambda^{\nu-1/2}\log
n,
\]
which diverges as $n\rightarrow\infty$. Constructions adapted to the fixed
terminal point are therefore required. In the empirical case this is achieved
through the conditional decomposition at $t$, whereas in the quantile case it
is achieved through the R\'enyi representation and the reversal of the
corresponding exponential spacings.

\smallskip\noindent To obtain the empirical-process construction, we split the
original uniform sample at $t$. The indicators
\[
\varepsilon_{i}:=\mathbb{I}\{U_{i}\leq t\},\qquad i\geq1,
\]
form an i.i.d.\ Bernoulli sequence. Moreover, the observations falling to the
left and to the right of $t$, after rescaling and retaining their respective
orders of appearance, form two independent i.i.d.\ uniform samples that are
also independent of the Bernoulli marks. This decomposition yields exact
representations of both the ordinary empirical process and its increments in
terms of two conditional uniform empirical processes and the common binomial
fluctuation at $t$.

\smallskip

\noindent The two conditional empirical processes are coupled with independent
Brownian bridges by means of the weighted approximation of \cite{CsCsHM86},
whereas the centered Bernoulli partial sums are coupled with a Wiener process
through the Koml\'{o}s--Major--Tusn\'{a}dy approximation. The three Gaussian
components are then pasted together at $t$. The resulting process is a
standard Brownian bridge and simultaneously provides the ordinary
empirical-process approximation and the fixed-terminal-point
empirical-increment approximation stated in Theorem~\ref{Theo1}.

\smallskip

\noindent The quantile-increment approximation is obtained through a parallel,
but structurally different, argument. By the R\'enyi representation, spacings
between uniform order statistics can be expressed as normalized sums of
independent standard exponential random variables. Reversing the exponential
spacings associated with the fixed terminal point $t$ transforms the quantile
increment process, on the empirical-quantile grid, into an ordinary uniform
quantile process. The weighted quantile-process approximation of
\cite{CsCsHM86} may then be applied down to the finite-count scale $\lambda
/n$. The remaining discretization and interpolation terms are controlled with
the same weight.

\smallskip

\noindent A preliminary unpublished version of the empirical-increment result
was announced by Necir (2017, arXiv:1709.00747). The present paper provides a
new proof starting from an arbitrary uniform sample, possibly after enlarging
its underlying probability space, and supplements it with the corresponding
weighted approximation for the uniform quantile increment process. The two
proofs reveal complementary mechanisms: rotation and conditional decomposition
for the empirical increments, and the R\'enyi representation with reversal of
exponential spacings for the quantile increments.

\smallskip

\noindent In Section~\ref{Sec2}, the empirical part of the main result is
applied to empirical subdistribution-tail processes arising from randomly
right-censored observations. The resulting approximations provide
probabilistic tools for the subsequent analysis of censored tail processes and
related extreme-value procedures. We now state the main theorem.

\begin{theorem}
\label{Theo1} Let $U_{1},U_{2},\ldots$ be any sequence of i.i.d.\ random
variables uniformly distributed on $(0,1)$ and defined on a probability space
$(\Omega_{0},\mathcal{A}_{0},\mathbb{P}_{0})$. Fix $t\in(0,1)$. There exists
an extension $(\Omega,\mathcal{A},\mathbb{P})$ of this probability space
carrying the original sequence $U_{1},U_{2},\ldots$, together with a sequence
of common standard Brownian bridges
\[
\left\{  B_{n}(u);\,0\leq u\leq1\right\}  _{n\geq1},
\]

such that, for every fixed $\lambda>0$ and every $0\leq\nu<1/4$,
\begin{equation}
\sup_{\lambda/n\leq u\leq1-\lambda/n}\frac{n^{\nu}\left\vert \alpha
_{n}(u)-B_{n}(u)\right\vert }{\left\{  u(1-u)\right\}  ^{1/2-\nu}%
}=O_{\mathbb{P}}(1),\qquad n\rightarrow\infty, \label{ordinary-approximation}%
\end{equation}
and
\begin{equation}
\sup_{\lambda/n\leq s<t}\frac{n^{\nu}\left\vert \alpha_{n}(s;t)-B_{n}%
(s;t)\right\vert }{s^{1/2-\nu}}=O_{\mathbb{P}}(1),\qquad n\rightarrow\infty,
\label{increment-approximation}%
\end{equation}
where
\[
B_{n}(s;t):=B_{n}(t)-B_{n}(t-s).
\]

Moreover, for every fixed $\lambda>0$ and every $0\leq\eta<1/2$,
\begin{equation}
\sup_{\lambda/n\leq u\leq1-\lambda/n}\frac{n^{\eta}\left\vert \beta
_{n}(u)-B_{n}(u)\right\vert }{\left\{  u(1-u)\right\}  ^{1/2-\eta}%
}=O_{\mathbb{P}}(1),\qquad n\rightarrow\infty,
\label{ordinary-quantile-approximation}%
\end{equation}
and
\begin{equation}
\sup_{\lambda/n\leq s<t}\frac{n^{\eta}\left\vert \beta_{n}(s;t)-B_{n}%
(s;t)\right\vert }{s^{1/2-\eta}}=O_{\mathbb{P}}(1),\qquad n\rightarrow\infty,
\label{quantile-increment-approximation}%
\end{equation}
In the last assertion, the notation $B_{n}(s;t)$ is used for a Gaussian
approximating process having the same distribution as the Brownian-bridge
increment $\{B_{n}(t)-B_{n}(t-s):0\leq s\leq t\}$; the auxiliary coupling is
made explicit in the proof. The extension of the probability space and the
common Brownian bridge may be chosen independently of the particular values of
$\lambda$, $\nu$, and $\eta$.
\end{theorem}

\begin{remark}
The increment approximation \eqref{increment-approximation} strengthens the
fixed-terminal-point approximation \eqref{Alex-emp-incr} by providing an
explicit stochastic rate. Indeed, let $0\leq\rho<1/2$, and choose
\[
0<\nu<\min\left\{  \frac14,\frac12-\rho\right\}  .
\]
Set
\[
T_{n}^{\left(  e\right)  } := \sup_{\lambda/n\leq s<t} \frac{n^{\nu}\left\vert
\alpha_{n}(s;t)-B_{n}(s;t) \right\vert } {s^{1/2-\nu}}.
\]
Then $T_{n}^{\left(  e\right)  }=O_{\mathbb{P}}(1)$, and
\[
\sup_{\lambda/n\leq s<t} \frac{\left\vert \alpha_{n}(s;t)-B_{n}%
(s;t)\right\vert } {s^{\rho}} \leq n^{-\nu}t^{1/2-\nu-\rho}T_{n}^{\left(
e\right)  } = o_{\mathbb{P}}(1).
\]
Thus, \eqref{increment-approximation} implies the corresponding $s^{-\rho}%
$-weighted convergence for every $\rho<1/2$. In addition, the case $\nu=0$
yields the critical bound
\[
\sup_{\lambda/n\leq s<t} \frac{\left\vert \alpha_{n}(s;t)-B_{n}%
(s;t)\right\vert } {s^{1/2}} = O_{\mathbb{P}}(1).
\]

\end{remark}

\begin{remark}
The quantile-increment approximation \eqref{quantile-increment-approximation}
is the fixed-terminal-point counterpart of the global increment approximation
of \cite[Theorem~4.6.1]{CsCsHM86}. Whereas the latter is formulated above the
logarithmic scale $c(\log n)/n$, assertion
\eqref{quantile-increment-approximation} reaches the finer finite-count scale
$\lambda/n$. In particular, if $0\leq\rho<1/2$ and
\[
0<\eta<\frac12-\rho,
\]
then
\[
\sup_{\lambda/n\leq s<t} \frac{\left\vert \beta_{n}(s;t)-B_{n}(s;t)
\right\vert } {s^{\rho}} = o_{\mathbb{P}}(1).
\]

\end{remark}

\begin{remark}
The restriction $0\leq\nu<1/4$ in
\eqref{ordinary-approximation}--\eqref{increment-approximation} is inherited
from the weighted empirical-process approximation of \cite{CsCsHM86} used in
the proof. At the boundary $\nu=1/4$, the corresponding source result contains
an additional logarithmic factor. The range $0\leq\eta<1/2$ in
\eqref{ordinary-quantile-approximation}--\eqref{quantile-increment-approximation}
is the classical range for the weighted quantile-process approximation. No
claim of optimality for these restrictions is made here. The constants
implicit in the $O_{\mathbb{P}}(1)$ terms may depend on the fixed parameters
$t$, $\lambda$, $\nu$, and $\eta$.
\end{remark}

\begin{remark}
The same sequence $\{B_{n}\}_{n\geq1}$ is used for the ordinary empirical
process, the ordinary quantile process, and the fixed-endpoint empirical
increments. This common Brownian-bridge notation is essential for keeping the
Gaussian dependence structure transparent in the censoring application. For
the quantile-increment assertion, the proof first constructs an auxiliary
Brownian bridge attached to the reversed quantile process and then identifies
its law with that of the increment process $\{B_{n}(t)-B_{n}(t-s):0\leq s\leq
t\}$. Thus the theorem is written in terms of the same Brownian-increment law,
while the proof keeps track of the auxiliary coupling needed to justify this representation.
\end{remark}

\section{Application to Censored Extreme-Value Statistics}

\label{Sec2}

\noindent We now illustrate the usefulness of Theorem~\ref{Theo1} in the
random right-censoring framework. We first derive exact representations of the
two empirical subdistribution-tail processes in terms of a single uniform
empirical process. We then apply, respectively, the ordinary empirical-process
and empirical-increment approximations of Theorem~\ref{Theo1}. This yields
approximations driven by a common Brownian bridge for the two subdistribution
components and therefore preserves their joint dependence structure.\smallskip

\noindent Let $X_{1},\ldots,X_{n}$ be independent copies of a nonnegative
random variable $X$ with continuous distribution function $F$, and let
$Y_{1},\ldots,Y_{n}$ be independent copies of a nonnegative random variable
$Y$ with continuous distribution function $G$. Assume that the two samples are
independent. For $1\leq j\leq n$, we observe
\[
Z_{j}:=\min(X_{j},Y_{j})
\]
together with the censoring indicator
\[
\delta_{j}:=\mathbb{I}\{X_{j}\leq Y_{j}\}.
\]
Thus, $\delta_{j}=1$ means that $X_{j}$ is observed, whereas $\delta_{j}=0$
means that it is right censored. If $H$ denotes the distribution function of
\[
Z:=\min(X,Y),
\]
then, by independence,
\[
1-H(x)=\{1-F(x)\}\{1-G(x)\},\qquad x\geq0.
\]
Define the subdistribution functions
\[
H^{(i)}(z):=\mathbb{P}\{Z\leq z,\delta=i\},\qquad z\geq0,\quad i=0,1.
\]
Then
\[
H(z)=H^{(0)}(z)+H^{(1)}(z).
\]
Their empirical counterparts are
\[
H_{n}^{(0)}(z):=\frac{1}{n}\sum_{j=1}^{n} \mathbb{I}\{Z_{j}\leq z\}(1-\delta
_{j})
\]
and
\[
H_{n}^{(1)}(z):=\frac{1}{n}\sum_{j=1}^{n} \mathbb{I}\{Z_{j}\leq z\}\delta
_{j}.
\]
Put
\[
\theta:=H^{(1)}(\infty)=\mathbb{P}\{\delta=1\},
\]
and assume throughout that
\[
0<\theta<1.
\]
Consequently,
\[
H^{(0)}(\infty)=1-\theta.
\]

\noindent We next use the uniform representation of \cite{EnKo92}. Define
\[
\xi_{j}:=\delta_{j}H^{(1)}(Z_{j})+(1-\delta_{j}) \left\{  \theta+H^{(0)}%
(Z_{j})\right\}  , \qquad j=1,\ldots,n.
\]
Conditionally on $\delta_{j}=1$, the distribution function of $Z_{j}$ is the
normalized subdistribution function
\[
z\longmapsto\mathbb{P}\{Z_{j}\leq z\mid\delta_{j}=1\} = \frac{H^{(1)}%
(z)}{\theta}.
\]
Since this conditional distribution function is continuous,
\[
\frac{H^{(1)}(Z_{j})}{\theta}
\]
is uniformly distributed on $(0,1)$, conditionally on $\delta_{j}=1$.
Equivalently, $H^{(1)}(Z_{j})$ is conditionally uniform on $(0,\theta)$.

\smallskip

\noindent Similarly, conditionally on $\delta_{j}=0$, the distribution
function of $Z_{j}$ is
\[
z\longmapsto\mathbb{P}\{Z_{j}\leq z\mid\delta_{j}=0\} = \frac{H^{(0)}%
(z)}{1-\theta},
\]
and hence
\[
\frac{H^{(0)}(Z_{j})}{1-\theta}
\]
is uniformly distributed on $(0,1)$. Therefore, $\theta+H^{(0)}(Z_{j})$ is
conditionally uniform on $(\theta,1)$.

\smallskip

\noindent It follows that, for $0\leq u\leq\theta$,
\[
\mathbb{P}\{\xi_{j}\leq u\} = \theta\frac{u}{\theta} = u,
\]
whereas, for $\theta<u\leq1$,
\[
\mathbb{P}\{\xi_{j}\leq u\} = \theta+(1-\theta)\frac{u-\theta}{1-\theta} = u.
\]
Consequently,
\[
\xi_{1},\ldots,\xi_{n}
\]
are i.i.d.\ random variables uniformly distributed on $(0,1)$.

\smallskip

\noindent Let
\[
\mathbb{U}_{n}(u) := \frac{1}{n}\sum_{j=1}^{n}\mathbb{I}\{\xi_{j}\leq u\},
\qquad0\leq u\leq1,
\]
and define the associated uniform empirical process by
\[
\alpha_{n}^{\ast}(u) := \sqrt{n}\{\mathbb{U}_{n}(u)-u\}, \qquad0\leq u\leq1.
\]

\smallskip

\noindent Following \cite{DeEn96}, the two empirical subdistribution-tail
processes admit the following almost sure representations in terms of the
uniform empirical process $\alpha_{n}^{\ast}$. For $i=0,1$, write
\[
\overline{H}^{(i)}(v) := H^{(i)}(\infty)-H^{(i)}(v)
\]
and
\[
\overline{H}_{n}^{(i)}(v) := H_{n}^{(i)}(\infty)-H_{n}^{(i)}(v), \qquad
v\geq0.
\]
Since
\[
1-\overline{H}^{(0)}(v) = \theta+H^{(0)}(v),
\]
the definition of $\xi_{j}$ gives, almost surely,
\[
\mathbb{U}_{n}\left(  1-\overline{H}^{(0)}(v)\right)  = 1-\overline{H}%
_{n}^{(0)}(v).
\]
Therefore,
\begin{equation}
\alpha_{n}^{\ast}\left(  1-\overline{H}^{(0)}(v)\right)  = -\sqrt{n}\left\{
\overline{H}_{n}^{(0)}(v)-\overline{H}^{(0)}(v) \right\}  .
\label{censored-tail-representation}%
\end{equation}
Similarly,
\[
\theta-\overline{H}^{(1)}(v) = H^{(1)}(v).
\]
Moreover, almost surely,
\[
\mathbb{U}_{n}(\theta) - \mathbb{U}_{n}\left(  \theta-\overline{H}%
^{(1)}(v)\right)  = \overline{H}_{n}^{(1)}(v).
\]
It follows that
\begin{equation}
\alpha_{n}^{\ast}\left(  \overline{H}^{(1)}(v);\theta\right)  = \sqrt
{n}\left\{  \overline{H}_{n}^{(1)}(v)-\overline{H}^{(1)}(v) \right\}  ,
\label{uncensored-tail-representation}%
\end{equation}
where
\[
\alpha_{n}^{\ast}(s;\theta) := \alpha_{n}^{\ast}(\theta)-\alpha_{n}^{\ast
}(\theta-s), \qquad0\leq s\leq\theta.
\]

\noindent For notational convenience, define
\[
\mathbb{H}_{n}^{(i)}(v) := \sqrt{n}\left\{  \overline{H}_{n}^{(i)}%
(v)-\overline{H}^{(i)}(v) \right\}  , \qquad i=0,1.
\]
Then \eqref{censored-tail-representation} and
\eqref{uncensored-tail-representation} become
\begin{equation}
\mathbb{H}_{n}^{(0)}(v) = -\alpha_{n}^{\ast}\left(  1-\overline{H}%
^{(0)}(v)\right)  \label{H0-representation}%
\end{equation}
and
\begin{equation}
\mathbb{H}_{n}^{(1)}(v) = \alpha_{n}^{\ast}\left(  \overline{H}^{(1)}%
(v);\theta\right)  . \label{H1-representation}%
\end{equation}
We now apply Theorem~\ref{Theo1} to the uniform empirical process $\alpha
_{n}^{\ast}$, with terminal point $t=\theta$. Possibly after enlarging the
underlying probability space, there exists a sequence of standard Brownian
bridges
\[
\left\{  \mathcal{B}_{n}(u);\,0\leq u\leq1\right\}  _{n\geq1}
\]
such that the ordinary empirical-process approximation and the
empirical-increment approximation below hold with the same bridge. The
possible dependence of $\mathcal{B}_{n}$ on the fixed terminal point $\theta$
is suppressed from the notation.

\smallskip

\noindent For every fixed $\lambda>0$ and every $0\leq\nu<1/4$, the interval
\[
\left[  \theta,1-\frac{\lambda}{n}\right]
\]
is contained in
\[
\left[  \frac{\lambda}{n},1-\frac{\lambda}{n}\right]
\]
for all sufficiently large $n$. Hence, restricting the ordinary empirical-
process approximation of Theorem~\ref{Theo1} gives
\begin{equation}
\sup_{\theta\leq u\leq1-\lambda/n} \frac{n^{\nu}\left\vert \alpha_{n}^{\ast
}(u)-\mathcal{B}_{n}(u) \right\vert } {(1-u)^{1/2-\nu}} = O_{\mathbb{P}}(1).
\label{ordinary-censored-approximation}%
\end{equation}
Indeed, on this interval,
\[
\frac{\left\vert \alpha_{n}^{\ast}(u)-\mathcal{B}_{n}(u) \right\vert }
{(1-u)^{1/2-\nu}} = u^{1/2-\nu} \frac{\left\vert \alpha_{n}^{\ast
}(u)-\mathcal{B}_{n}(u) \right\vert } {\left\{  u(1-u)\right\}  ^{1/2-\nu}},
\]
and $u^{1/2-\nu}\leq1$.

\smallskip

\noindent The empirical-increment part of Theorem~\ref{Theo1} yields
\begin{equation}
\sup_{\lambda/n\leq s<\theta} \frac{n^{\nu}\left\vert \alpha_{n}^{\ast
}(s;\theta)-\mathcal{B}_{n}(s;\theta) \right\vert } {s^{1/2-\nu}} =
O_{\mathbb{P}}(1), \label{increment-censored-approximation}%
\end{equation}
where
\[
\mathcal{B}_{n}(s;\theta) := \mathcal{B}_{n}(\theta)-\mathcal{B}_{n}%
(\theta-s).
\]
We first apply \eqref{ordinary-censored-approximation} to the censored
component. Taking
\[
u=1-\overline{H}^{(0)}(v)
\]
and using \eqref{H0-representation}, we obtain
\begin{equation}
\sup_{\substack{v\geq0:\\\lambda/n\leq\overline{H}^{(0)}(v)\leq1-\theta}}
\frac{n^{\nu}\left\vert \mathbb{H}_{n}^{(0)}(v) + \mathcal{B}_{n}\left(
1-\overline{H}^{(0)}(v)\right)  \right\vert } {\left\{  \overline{H}%
^{(0)}(v)\right\}  ^{1/2-\nu}} = O_{\mathbb{P}}(1).
\label{censored-tail-approximation}%
\end{equation}
Next, taking
\[
s=\overline{H}^{(1)}(v)
\]
in \eqref{increment-censored-approximation} and using
\eqref{H1-representation}, we obtain
\begin{equation}
\sup_{\substack{v\geq0:\\\lambda/n\leq\overline{H}^{(1)}(v)<\theta}}
\frac{n^{\nu}\left\vert \mathbb{H}_{n}^{(1)}(v) - \mathcal{B}_{n}\left(
\overline{H}^{(1)}(v);\theta\right)  \right\vert } {\left\{  \overline
{H}^{(1)}(v)\right\}  ^{1/2-\nu}} = O_{\mathbb{P}}(1).
\label{uncensored-tail-approximation}%
\end{equation}
Relations \eqref{censored-tail-approximation} and
\eqref{uncensored-tail-approximation} show that the two empirical
subdistribution-tail processes are approximated by functionals of the same
Brownian bridge. More precisely, the censored component $\mathbb{H}_{n}%
^{(0)}(v)$ is approximated by
\[
-\mathcal{B}_{n}\left(  1-\overline{H}^{(0)}(v)\right)  ,
\]
whereas the uncensored component $\mathbb{H}_{n}^{(1)}(v)$ is approximated by
the bridge increment
\[
\mathcal{B}_{n}\left(  \overline{H}^{(1)}(v);\theta\right)  .
\]
The common coupling therefore preserves the joint Gaussian dependence between
the two components.

\smallskip

\noindent These approximations provide process-level ingredients for the
asymptotic analysis of statistics involving both empirical subdistribution
tails. Combined with suitable functional representations and appropriate
controls of bias and remainder terms, they can be used to derive limit
distributions for censored tail processes, tail product-limit procedures, and
related extreme-value estimators; see, for example,
\cite{BMN2015,MNS2025,GNM2026A,GNM2026B}.

\section{Proof of Theorem~\ref{Theo1}\label{sec3}}

\noindent Fix $t\in(0,1)$. We first construct, possibly on an extension of the
original probability space, the common Brownian bridge used for the ordinary
empirical process, the ordinary quantile process, and the empirical
increments. For the quantile increments, an auxiliary bridge is introduced in
the proof and then identified in distribution with the increment of a standard
Brownian bridge. These constructions may depend on the fixed terminal point
$t$, but they do not depend on the parameters $\lambda$, $\nu$, and $\eta
$.\smallskip

\noindent Once the constructions have been completed, let $\lambda>0$,
$0\leq\nu<1/4$, and $0\leq\eta<1/2$ be fixed, and put
\[
\rho:=\frac{1}{2}-\nu\qquad\text{and}\qquad\kappa:=\frac{1}{2}-\eta.
\]
Throughout the proof, the possible dependence of the approximating Brownian
bridges on $t$ is suppressed from the notation whenever no confusion can
arise.\smallskip

\noindent The proof is divided into two parts. The first part establishes the
ordinary empirical-process approximation and the empirical-increment
approximation with the same sequence of Brownian bridges. The second part
proves the corresponding assertions for the uniform quantile process and its
increments.\smallskip

\noindent For the empirical part, we exploit the fixed terminal point $t$ by
splitting the original uniform sample according to whether each observation
falls in $(0,t]$ or in $(t,1)$. The indicators recording this dichotomy form
an i.i.d.\ Bernoulli sequence with success probability $t$, whose partial sums
count the observations lying in $(0,t]$. After rescaling and retaining their
respective orders of appearance, the observations falling to the left and to
the right of $t$ generate two independent i.i.d.\ uniform samples, which are
also independent of the Bernoulli marks.\smallskip

\noindent This splitting yields exact decompositions of both the ordinary
empirical process and its increments ending at $t$ into three independent
components: a left conditional empirical process, a right conditional
empirical process, and the centered Bernoulli fluctuation at $t$. We couple
the two conditional empirical processes with independent Brownian bridges by
means of appropriate one-sided weighted approximations, and couple the
centered Bernoulli partial sums with a Wiener process through the
Koml\'os--Major--Tusn\'ady approximation.\smallskip

\noindent The three Gaussian components are then pasted together at $t$ to
form a single standard Brownian bridge on $[0,1]$. The exact decompositions of
the ordinary empirical process and of its increments are compared with the
corresponding decompositions of this common bridge, and the resulting coupling
errors are controlled uniformly over the ranges appearing in
\eqref{ordinary-approximation} and \eqref{increment-approximation}.\smallskip

\noindent For the quantile part, the ordinary weighted approximation is taken
from the joint empirical--quantile construction of \cite{CsCsHM86}. With our
sign convention $\beta_{n}(u)=\sqrt n\{u-V_{n}(u)\}$, this construction uses
the same Brownian bridge $B_{n}$ as in the ordinary empirical approximation.
Thus no additional Brownian bridge is needed for the ordinary quantile
process. To treat the increments ending at $t$, we use the R\'enyi
representation of the uniform order statistics. After reversing the
exponential spacings associated with the fixed terminal point, the quantile
increment process is transformed, on the empirical-quantile grid, into an
ordinary uniform quantile process. The weighted approximation of
\cite{CsCsHM86} can then be applied down to the finite-count scale $\lambda
/n$, and the remaining discretization terms are controlled with the same
weight. This yields \eqref{ordinary-quantile-approximation} and \eqref{quantile-increment-approximation}.

\subsection{Uniform empirical process and its increments}

\mbox{}

\smallskip

\noindent\textbf{Step 1: Splitting the uniform sample at }$t$\textbf{.}
\smallskip

\noindent Let
\[
U_{1},U_{2},\ldots
\]
be the i.i.d.\ uniform random variables appearing in Theorem~\ref{Theo1},
defined on the initial probability space $(\Omega_{0},\mathcal{A}%
_{0},\mathbb{P}_{0})$. Define
\[
\varepsilon_{i}:=\mathbb{I}\{U_{i}\leq t\}, \qquad i\geq1.
\]
Then $\varepsilon_{1},\varepsilon_{2},\ldots$ are i.i.d.\ Bernoulli random
variables satisfying
\[
\mathbb{P}_{0}\{\varepsilon_{i}=1\}=t, \qquad\mathbb{P}_{0}\{\varepsilon
_{i}=0\}=1-t.
\]

\noindent Let
\[
\tau_{j}^{(L)} := \inf\left\{  k\geq1:\sum_{i=1}^{k}\varepsilon_{i}=j
\right\}  , \qquad j\geq1,
\]
and
\[
\tau_{j}^{(R)} := \inf\left\{  k\geq1:\sum_{i=1}^{k}(1-\varepsilon_{i})=j
\right\}  , \qquad j\geq1.
\]
Since $0<t<1$, the Bernoulli sequence contains infinitely many zeros and
infinitely many ones almost surely. Hence all these indices are almost surely
finite. Define the rescaled left and right subsequences by
\[
V_{j}^{(L)} := \frac{U_{\tau_{j}^{(L)}}}{t}, \qquad j\geq1,
\]
and
\[
V_{j}^{(R)} := \frac{U_{\tau_{j}^{(R)}}-t}{1-t}, \qquad j\geq1.
\]

\noindent Conditionally on the entire mark sequence $\{\varepsilon
_{i}\}_{i\geq1}$, the variables with mark $1$ are independent and uniformly
distributed on $(0,t)$, whereas those with mark $0$ are independent and
uniformly distributed on $(t,1)$. Consequently, for any positive integers $r$
and $q$, and any Borel subsets $A_{1},\ldots,A_{r},C_{1},\ldots,C_{q}$ of
$(0,1)$,
\begin{align*}
&  \mathbb{P}_{0}\left\{  V_{1}^{(L)}\in A_{1},\ldots,V_{r}^{(L)}\in A_{r},\,
V_{1}^{(R)}\in C_{1},\ldots,V_{q}^{(R)}\in C_{q} \,\big|\, \{\varepsilon
_{i}\}_{i\geq1} \right\} \\
&  \qquad= \prod_{j=1}^{r}\operatorname{Leb}(A_{j}) \prod_{k=1}^{q}%
\operatorname{Leb}(C_{k}), \qquad\text{a.s.}%
\end{align*}
The conditional distribution on the right-hand side does not depend on the
realized mark sequence. Therefore,
\[
\{\varepsilon_{i}\}_{i\geq1}, \qquad\{V_{j}^{(L)}\}_{j\geq1}, \qquad
\{V_{j}^{(R)}\}_{j\geq1}
\]
are mutually independent random elements, and the two $V$-sequences are
independent i.i.d.\ uniform samples on $(0,1)$.

\noindent For $m\geq1$, define
\[
\mathbb{V}_{m}^{(L)}(x) := \frac{1}{m}\sum_{j=1}^{m} \mathbb{I}\{V_{j}%
^{(L)}\leq x\}, \qquad0\leq x\leq1,
\]
and
\[
\mathbb{V}_{m}^{(R)}(x) := \frac{1}{m}\sum_{j=1}^{m} \mathbb{I}\{V_{j}%
^{(R)}\leq x\}, \qquad0\leq x\leq1.
\]
Their associated uniform empirical processes are
\[
\gamma_{m}^{(L)}(x) := \sqrt m\left\{  \mathbb{V}_{m}^{(L)}(x)-x \right\}
\]
and
\[
\gamma_{m}^{(R)}(x) := \sqrt m\left\{  \mathbb{V}_{m}^{(R)}(x)-x \right\}  .
\]
We adopt the conventions
\[
\mathbb{V}_{0}^{(L)} = \mathbb{V}_{0}^{(R)} \equiv0, \qquad\gamma_{0}^{(L)} =
\gamma_{0}^{(R)} \equiv0.
\]

\noindent\textbf{Step 2: Coupling the three independent components.}
\smallskip

\noindent We first attach Brownian bridges to the two conditional uniform
samples. Corollary~2.1 of \cite{CsCsHM86} provides a coupling of an
i.i.d.\ uniform sequence and a sequence of standard Brownian bridges
satisfying the weighted approximations used below.

\noindent The spaces
\[
(0,1)^{\mathbb{N}} \qquad\text{and}\qquad C[0,1]^{\mathbb{N}}
\]
are standard Borel spaces. Hence the joint distribution supplied by
\cite{CsCsHM86} admits a regular conditional distribution of the sequence of
bridges given the uniform sequence. By the disintegration and randomization
theorem for probability kernels, after enlarging the original probability
space, we may attach to $\{V_{j}^{(L)}\}_{j\geq1}$ a sequence of standard
Brownian bridges
\[
\left\{  L_{m}(x);\,0\leq x\leq1 \right\}  _{m\geq1}
\]
such that
\[
\left(  \{V_{j}^{(L)}\}_{j\geq1}, \{L_{m}\}_{m\geq1} \right)
\]
has the same joint distribution as the coupling constructed in \cite{CsCsHM86}%
. Using an independent auxiliary randomization, we may similarly attach to
$\{V_{j}^{(R)}\}_{j\geq1}$ a sequence
\[
\left\{  R_{m}(x);\,0\leq x\leq1 \right\}  _{m\geq1}.
\]
The resulting left and right empirical-process couplings are mutually
independent and are independent of the Bernoulli sequence $\{\varepsilon
_{i}\}_{i\geq1}$.

\noindent For $m\geq1$, put
\[
\Delta_{m}^{(L)}(x) := \gamma_{m}^{(L)}(x)-L_{m}(x)
\]
and
\[
\Delta_{m}^{(R)}(x) := \gamma_{m}^{(R)}(x)-R_{m}(x).
\]
By relations~(2.30) and~(2.31) of Corollary~2.1 in \cite{CsCsHM86}, for every
fixed $a>0$, every $0\leq\nu<1/4$, and each $J\in\{L,R\}$,
\begin{equation}
\sup_{a/m\leq x\leq1} \frac{m^{\nu}\left\vert \Delta_{m}^{(J)}(x)\right\vert }
{x^{\rho}} = O_{\mathbb{P}}(1), \qquad m\rightarrow\infty,
\label{near-zero-weighted-approximation}%
\end{equation}
and
\begin{equation}
\sup_{0\leq x\leq1-a/m} \frac{m^{\nu}\left\vert \Delta_{m}^{(J)}(x)\right\vert
} {(1-x)^{\rho}} = O_{\mathbb{P}}(1), \qquad m\rightarrow\infty.
\label{near-one-weighted-approximation}%
\end{equation}
The underlying couplings may be fixed independently of the particular values
of $a$ and $\nu$.

\noindent Let $L_{0}$ and $R_{0}$ be two additional independent standard
Brownian bridges, independent of all the preceding random elements.

We next couple the Bernoulli fluctuation. Put
\[
S_{n} := \sum_{i=1}^{n}(\varepsilon_{i}-t), \qquad\sigma_{t}^{2}:=t(1-t).
\]
\cite{KMT75,KMT76} supply a joint distribution of an i.i.d.\ Bernoulli
sequence with success probability $t$ and a standard Wiener process satisfying
the corresponding strong approximation. Applying the same disintegration
argument, on a further extension if necessary, we may attach to the already
given Bernoulli sequence a standard Wiener process $W$. The auxiliary
randomization may be chosen independently of those used for the left and right
empirical-process couplings. Consequently, the pair
\[
\left(  \{\varepsilon_{i}\}_{i\geq1},W \right)
\]
is independent of the two empirical-process couplings, and
\begin{equation}
\max_{1\leq k\leq n} \left\vert S_{k}-\sigma_{t}W(k) \right\vert = O(\log n),
\qquad\text{a.s.} \label{bernoulli-kmt}%
\end{equation}

\noindent Define
\[
N_{n} := \sum_{i=1}^{n}\varepsilon_{i} = \sum_{i=1}^{n}\mathbb{I}\{U_{i}\leq
t\}.
\]
Then
\[
\frac{N_{n}}{n} \longrightarrow t, \qquad\text{a.s.},
\]
and
\[
Z_{n} := \frac{N_{n}-nt}{\sqrt n} = \alpha_{n}(t).
\]
Set
\[
G_{n} := \frac{\sigma_{t}W(n)}{\sqrt n}.
\]
Then
\[
G_{n} \sim\mathcal{N}\bigl(0,t(1-t)\bigr),
\]
and \eqref{bernoulli-kmt} implies
\begin{equation}
\left\vert Z_{n}-G_{n}\right\vert = O\left(  \frac{\log n}{\sqrt n}\right)  ,
\qquad\text{a.s.} \label{binomial-normal-approximation}%
\end{equation}
In particular,
\[
Z_{n} = O_{\mathbb{P}}(1).
\]

\noindent We shall use the following elementary random-index fact. Suppose
that
\[
X_{m} = O_{\mathbb{P}}(1), \qquad m\rightarrow\infty,
\]
that $M_{n}\rightarrow\infty$ in probability, and that $M_{n}$ is independent
of the entire sequence $\{X_{m}\}_{m\geq1}$. Then
\[
X_{M_{n}} = O_{\mathbb{P}}(1).
\]
Indeed, for every $\epsilon>0$, there exist $C<\infty$ and $m_{0}\geq1$ such
that
\[
\sup_{m\geq m_{0}} \mathbb{P}\left\{  \left\vert X_{m}\right\vert >C\right\}
< \epsilon.
\]
Therefore,
\[
\mathbb{P}\left\{  \left\vert X_{M_{n}}\right\vert >C\right\}  \leq
\mathbb{P}\{M_{n}<m_{0}\} + \sup_{m\geq m_{0}} \mathbb{P}\left\{  \left\vert
X_{m}\right\vert >C\right\}  ,
\]
and the right-hand side is eventually smaller than $2\epsilon$.

\noindent Since
\[
N_{n}\longrightarrow\infty\qquad\text{and}\qquad n-N_{n}\longrightarrow\infty
\]
in probability, and since these indices are independent of the left and right
empirical-process couplings, respectively, the preceding observation allows us
to use \eqref{near-zero-weighted-approximation} and
\eqref{near-one-weighted-approximation} with the random indices $N_{n}$ and
$n-N_{n}$.

\noindent For convenience, write
\[
L_{n}^{\ast} := L_{N_{n}}, \qquad R_{n}^{\ast} := R_{n-N_{n}}.
\]
We now verify the joint distribution of these randomly indexed bridges. Let
$A,D$ be Borel subsets of $C[0,1]$, and let $C_{0}$ be an event in
$\sigma(N_{n},G_{n})$. The two bridge families are mutually independent and
independent of $\sigma(N_{n},G_{n})$. Moreover, the distribution of $L_{m}$
and that of $R_{m}$ do not depend on $m$. Hence
\begin{align*}
&  \mathbb{P}\left\{  L_{N_{n}}\in A,\, R_{n-N_{n}}\in D,\, C_{0} \right\} \\
&  \qquad= \sum_{m=0}^{n} \mathbb{P}\{L_{m}\in A\} \mathbb{P}\{R_{n-m}\in D\}
\mathbb{P}\{C_{0},N_{n}=m\}\\
&  \qquad= \mathbb{P}\{L_{0}\in A\} \mathbb{P}\{R_{0}\in D\} \mathbb{P}%
(C_{0}).
\end{align*}
Thus, $L_{n}^{\ast}$ and $R_{n}^{\ast}$ are independent standard Brownian
bridges, and they are independent of $G_{n}$.

\noindent Finally, among $U_{1},\ldots,U_{n}$, the observations falling in
$(0,t]$, viewed as an unordered collection, are precisely
\[
tV_{1}^{(L)},\ldots,tV_{N_{n}}^{(L)},
\]
whereas those falling in $(t,1)$ are precisely
\[
t+(1-t)V_{1}^{(R)},\ldots, t+(1-t)V_{n-N_{n}}^{(R)}.
\]

\smallskip

\noindent\textbf{Step 3: Exact decompositions of the empirical process.}
\smallskip

\noindent The preceding splitting gives exact representations of the empirical
process on both sides of $t$.

\noindent Let $0\leq u\leq t$, and put
\[
x:=\frac{u}{t}.
\]
Then
\[
U_{n}(u) = \frac{N_{n}}{n}\mathbb{V}_{N_{n}}^{(L)}(x),
\]
where the equality remains valid when $N_{n}=0$, according to the convention
introduced above. Therefore,
\begin{align}
\alpha_{n}(u)  &  = \sqrt{n}\left\{  \frac{N_{n}}{n}\mathbb{V}_{N_{n}}%
^{(L)}(x)-tx \right\} \nonumber\\
&  = \sqrt{\frac{N_{n}}{n}}\, \gamma_{N_{n}}^{(L)}(x) + xZ_{n}.
\label{left-empirical-decomposition}%
\end{align}

Now let $t\leq u\leq1$, and put
\[
y:=\frac{u-t}{1-t}.
\]
Then
\[
U_{n}(u) = \frac{N_{n}}{n} + \frac{n-N_{n}}{n} \mathbb{V}_{n-N_{n}}^{(R)}(y).
\]
Since
\[
u=t+(1-t)y,
\]
we obtain
\begin{equation}
\alpha_{n}(u) = \sqrt{\frac{n-N_{n}}{n}}\, \gamma_{n-N_{n}}^{(R)}(y) +
(1-y)Z_{n}. \label{right-empirical-decomposition}%
\end{equation}

\noindent For $0\leq s\leq t$, put
\[
x:=\frac{s}{t}.
\]
Since $t-s=t(1-x)$, relation \eqref{left-empirical-decomposition} gives
\[
\alpha_{n}(t-s) = \sqrt{\frac{N_{n}}{n}}\, \gamma_{N_{n}}^{(L)}(1-x) +
(1-x)Z_{n}.
\]
Because $\alpha_{n}(t)=Z_{n}$, it follows that
\begin{equation}
\alpha_{n}(s;t) = xZ_{n} - \sqrt{\frac{N_{n}}{n}}\, \gamma_{N_{n}}^{(L)}(1-x),
\qquad x=\frac{s}{t}. \label{increment-exact-decomposition}%
\end{equation}

\smallskip\noindent\textbf{Step 4: Construction of the common Brownian
bridge.} \smallskip

\noindent We now paste the three Gaussian components together in a manner that
mirrors the preceding empirical decompositions. Define, for $0\leq u\leq1$,
\[
B_{n}(u):=
\begin{cases}
\displaystyle \frac{u}{t}G_{n} +\sqrt{t}\,L_{n}^{\ast}\left(  \frac{u}%
{t}\right)  , & 0\leq u\leq t,\\[1.2em]%
\displaystyle \frac{1-u}{1-t}G_{n} +\sqrt{1-t}\,R_{n}^{\ast}\left(  \frac
{u-t}{1-t}\right)  , & t<u\leq1.
\end{cases}
\]
At $u=t$, both expressions equal $G_{n}$, because
\[
L_{n}^{\ast}(1)=R_{n}^{\ast}(0)=0.
\]
Therefore, $B_{n}$ has continuous sample paths. Moreover,
\[
B_{n}(0)=B_{n}(1)=0.
\]

\noindent Since $G_{n}$, $L_{n}^{\ast}$, and $R_{n}^{\ast}$ are independent
centered Gaussian random elements, $B_{n}$ is a centered Gaussian process. We
identify its covariance function.

\noindent Let $0\leq u\leq v\leq t$, and write
\[
x:=\frac{u}{t}, \qquad z:=\frac{v}{t}.
\]
Using $x\leq z$, we obtain
\begin{align*}
\operatorname{Cov}\left\{  B_{n}(u),B_{n}(v)\right\}   &  =
xz\,t(1-t)+t\left\{  x-xz\right\} \\
&  = u(1-v).
\end{align*}

Let $t\leq u\leq v\leq1$, and put
\[
x:=\frac{u-t}{1-t}, \qquad z:=\frac{v-t}{1-t}.
\]
Using $x\leq z$, we have
\begin{align*}
\operatorname{Cov}\left\{  B_{n}(u),B_{n}(v)\right\}   &  = (1-x)(1-z)t(1-t)
+(1-t)\left\{  x-xz\right\} \\
&  = (1-t)(1-z)\left\{  t+(1-t)x\right\} \\
&  = u(1-v).
\end{align*}

\noindent Finally, if $0\leq u\leq t<v\leq1$, the independence of the left and
right bridge components gives
\begin{align*}
\operatorname{Cov}\left\{  B_{n}(u),B_{n}(v)\right\}   &  = \frac{u}{t}%
\frac{1-v}{1-t}t(1-t)\\
&  = u(1-v).
\end{align*}
Consequently,
\[
\operatorname{Cov}\left\{  B_{n}(u),B_{n}(v)\right\}  = \min(u,v)-uv,
\qquad0\leq u,v\leq1.
\]
Thus, $B_{n}$ is a standard Brownian bridge.

\noindent The increment of this bridge ending at $t$ also admits a simple
representation. If $0\leq s\leq t$ and $x=s/t$, then
\[
B_{n}(t)=G_{n}
\]
and
\[
B_{n}(t-s) = (1-x)G_{n} +\sqrt{t}\,L_{n}^{\ast}(1-x).
\]
Hence
\begin{equation}
B_{n}(s;t) = B_{n}(t)-B_{n}(t-s) = xG_{n} -\sqrt{t}\,L_{n}^{\ast}(1-x).
\label{bridge-increment-decomposition}%
\end{equation}

\smallskip

\noindent\textbf{Step 5: Approximation of the ordinary empirical process.}
\smallskip

\noindent Define
\[
A_{n} := \left\{  \left\vert \frac{N_{n}}{n}-t\right\vert \leq\frac{1}{2}%
\min(t,1-t) \right\}  .
\]
By the strong law of large numbers,
\[
\mathbb{P}(A_{n}) \longrightarrow1.
\]
On $A_{n}$,
\begin{equation}
\frac{N_{n}}{n} \geq\frac{t}{2}, \qquad\frac{n-N_{n}}{n} \geq\frac{1-t}{2}.
\label{sample-size-lower-bounds}%
\end{equation}
Moreover, since
\[
\frac{N_{n}}{n}-t = \frac{Z_{n}}{\sqrt{n}}
\]
and $Z_{n}=O_{\mathbb{P}}(1)$,
\begin{equation}
\left\vert \sqrt{\frac{N_{n}}{n}}-\sqrt{t} \right\vert = O_{\mathbb{P}}\left(
n^{-1/2}\right)  \label{left-square-root-error}%
\end{equation}
and
\begin{equation}
\left\vert \sqrt{\frac{n-N_{n}}{n}}-\sqrt{1-t} \right\vert = O_{\mathbb{P}%
}\left(  n^{-1/2}\right)  . \label{right-square-root-error}%
\end{equation}
Indeed, on $A_{n}$, the denominators in the identities
\[
\left\vert \sqrt{\frac{N_{n}}{n}}-\sqrt{t} \right\vert = \frac{\left\vert
N_{n}/n-t\right\vert } {\sqrt{N_{n}/n}+\sqrt{t}}
\]
and
\[
\left\vert \sqrt{\frac{n-N_{n}}{n}}-\sqrt{1-t} \right\vert = \frac{\left\vert
N_{n}/n-t\right\vert } {\sqrt{(n-N_{n})/n}+\sqrt{1-t}}
\]
are bounded away from zero.

\noindent For all sufficiently large $n$,
\[
\frac{\lambda}{n} < \min(t,1-t).
\]
Put
\[
a:=\frac{\lambda}{2}.
\]

\noindent We first consider the interval to the left of $t$. Let
\[
\frac{\lambda}{n} \leq u \leq t, \qquad x:=\frac{u}{t}.
\]
On $A_{n}$, relation \eqref{sample-size-lower-bounds} gives
\[
N_{n}x = \frac{N_{n}u}{t} \geq\frac{\lambda N_{n}}{nt} \geq\frac{\lambda}{2} =
a.
\]
Thus,
\[
x \geq\frac{a}{N_{n}},
\]
and the random-index version of \eqref{near-zero-weighted-approximation} applies.

\noindent Combining \eqref{left-empirical-decomposition} with the definition
of $B_{n}$, we obtain
\begin{align}
\alpha_{n}(u)-B_{n}(u) ={}  &  \sqrt{\frac{N_{n}}{n}}\, \Delta_{N_{n}}%
^{(L)}(x)\nonumber\\
&  + \left(  \sqrt{\frac{N_{n}}{n}}-\sqrt{t} \right)  L_{n}^{\ast}(x) +
x(Z_{n}-G_{n}). \label{left-error-decomposition}%
\end{align}

\noindent For the first term, on $A_{n}$,
\begin{align*}
&  \sup_{\lambda/n\leq u\leq t} \frac{ n^{\nu}\sqrt{N_{n}/n}\, \left\vert
\Delta_{N_{n}}^{(L)}(u/t)\right\vert } {\left\{  u(1-u)\right\}  ^{\rho}}\\
&  \qquad\leq\left(  \frac{N_{n}}{n}\right)  ^{\rho} \sup_{a/N_{n}\leq x\leq1}
\frac{ N_{n}^{\nu} \left\vert \Delta_{N_{n}}^{(L)}(x)\right\vert } {x^{\rho}}
\sup_{0<x\leq1} \frac{x^{\rho}} {\left\{  tx(1-tx)\right\}  ^{\rho}}.
\end{align*}
\noindent Since
\[
1-tx \geq1-t,
\]
the last supremum is bounded by
\[
\left\{  t(1-t)\right\}  ^{-\rho}.
\]
The random-index observation and \eqref{near-zero-weighted-approximation}
therefore imply
\begin{equation}
\sup_{\lambda/n\leq u\leq t} \frac{ n^{\nu}\sqrt{N_{n}/n}\, \left\vert
\Delta_{N_{n}}^{(L)}(u/t)\right\vert } {\left\{  u(1-u)\right\}  ^{\rho}} =
O_{\mathbb{P}}(1). \label{left-leading-error}%
\end{equation}

Since $L_{n}^{\ast}$ is a standard Brownian bridge,
\[
\left\Vert L_{n}^{\ast}\right\Vert _{\infty} := \sup_{0\leq x\leq1} \left\vert
L_{n}^{\ast}(x)\right\vert = O_{\mathbb{P}}(1).
\]
Moreover,
\[
u(1-u) \geq(1-t)\frac{\lambda}{n}, \qquad\frac{\lambda}{n}\leq u\leq t.
\]
Using \eqref{left-square-root-error} and $\nu+\rho=1/2$, we obtain
\begin{equation}
\sup_{\lambda/n\leq u\leq t} \frac{ n^{\nu} \left\vert \sqrt{N_{n}/n}-\sqrt
{t}\right\vert \left\vert L_{n}^{\ast}(u/t)\right\vert } {\left\{
u(1-u)\right\}  ^{\rho}} = O_{\mathbb{P}}(1). \label{left-coefficient-error}%
\end{equation}

\noindent Furthermore,
\[
\frac{x}{\left\{  u(1-u)\right\}  ^{\rho}} = \frac{u/t}{\left\{
u(1-u)\right\}  ^{\rho}} = \frac{u^{1-\rho}}{t(1-u)^{\rho}} \leq C_{t},
\qquad\frac{\lambda}{n}\leq u\leq t.
\]
It follows from \eqref{binomial-normal-approximation} that
\begin{equation}
\sup_{\lambda/n\leq u\leq t} \frac{ n^{\nu}x\left\vert Z_{n}-G_{n}\right\vert
} {\left\{  u(1-u)\right\}  ^{\rho}} = O_{\mathbb{P}} \left(  n^{\nu-1/2}\log
n \right)  = o_{\mathbb{P}}(1). \label{left-binomial-error}%
\end{equation}

\noindent Since $\mathbb{P}(A_{n}^{c})\longrightarrow0$, relations
\eqref{left-error-decomposition}--\eqref{left-binomial-error} yield
\begin{equation}
\sup_{\lambda/n\leq u\leq t} \frac{ n^{\nu} \left\vert \alpha_{n}%
(u)-B_{n}(u)\right\vert } {\left\{  u(1-u)\right\}  ^{1/2-\nu}} =
O_{\mathbb{P}}(1). \label{left-ordinary-approximation}%
\end{equation}

\noindent We next consider the interval to the right of $t$. Let
\[
t \leq u \leq1-\frac{\lambda}{n}, \qquad y:=\frac{u-t}{1-t}.
\]
On $A_{n}$,
\[
(n-N_{n})(1-y) = \frac{(n-N_{n})(1-u)}{1-t} \geq\frac{\lambda(n-N_{n}%
)}{n(1-t)} \geq\frac{\lambda}{2} = a.
\]
Therefore,
\[
y \leq1-\frac{a}{n-N_{n}},
\]
and the random-index version of \eqref{near-one-weighted-approximation} applies.

\noindent By \eqref{right-empirical-decomposition},
\begin{align}
\alpha_{n}(u)-B_{n}(u) ={}  &  \sqrt{\frac{n-N_{n}}{n}}\, \Delta_{n-N_{n}%
}^{(R)}(y)\nonumber\\
&  + \left(  \sqrt{\frac{n-N_{n}}{n}}-\sqrt{1-t} \right)  R_{n}^{\ast
}(y)\nonumber\\
&  + (1-y)(Z_{n}-G_{n}). \label{right-error-decomposition}%
\end{align}

\noindent For the first term, on $A_{n}$,
\begin{align*}
&  \sup_{t\leq u\leq1-\lambda/n} \frac{ n^{\nu}\sqrt{(n-N_{n})/n}\, \left\vert
\Delta_{n-N_{n}}^{(R)}(y)\right\vert } {\left\{  u(1-u)\right\}  ^{\rho}}\\
&  \qquad\leq\left(  \frac{n-N_{n}}{n}\right)  ^{\rho} \sup_{0\leq
y\leq1-a/(n-N_{n})} \frac{ (n-N_{n})^{\nu} \left\vert \Delta_{n-N_{n}}%
^{(R)}(y)\right\vert } {(1-y)^{\rho}}\\
&  \hspace{4cm}\times\sup_{0\leq y<1} \frac{(1-y)^{\rho}} {\left[  \left\{
t+(1-t)y\right\}  (1-t)(1-y) \right]  ^{\rho}}.
\end{align*}
Since
\[
t+(1-t)y \geq t,
\]
the last supremum is bounded by
\[
\left\{  t(1-t)\right\}  ^{-\rho}.
\]
The random-index observation and \eqref{near-one-weighted-approximation}
imply
\begin{equation}
\sup_{t\leq u\leq1-\lambda/n} \frac{ n^{\nu}\sqrt{(n-N_{n})/n}\, \left\vert
\Delta_{n-N_{n}}^{(R)}(y)\right\vert } {\left\{  u(1-u)\right\}  ^{\rho}} =
O_{\mathbb{P}}(1). \label{right-leading-error}%
\end{equation}

\noindent Since $R_{n}^{\ast}$ is a standard Brownian bridge,
\[
\left\Vert R_{n}^{\ast}\right\Vert _{\infty} = O_{\mathbb{P}}(1).
\]
Moreover,
\[
u(1-u) \geq t\frac{\lambda}{n}, \qquad t\leq u\leq1-\frac{\lambda}{n}.
\]
Using \eqref{right-square-root-error}, we obtain
\begin{equation}
\sup_{t\leq u\leq1-\lambda/n} \frac{ n^{\nu} \left\vert \sqrt{(n-N_{n}%
)/n}-\sqrt{1-t} \right\vert \left\vert R_{n}^{\ast}(y)\right\vert } {\left\{
u(1-u)\right\}  ^{\rho}} = O_{\mathbb{P}}(1). \label{right-coefficient-error}%
\end{equation}

\noindent Finally,
\begin{align*}
\frac{1-y}{\left\{  u(1-u)\right\}  ^{\rho}}  &  = \frac{1-u} {(1-t)\left\{
u(1-u)\right\}  ^{\rho}}\\
&  = \frac{(1-u)^{1-\rho}} {(1-t)u^{\rho}} \leq C_{t}, \qquad t\leq
u\leq1-\frac{\lambda}{n}.
\end{align*}
Therefore,
\begin{equation}
\sup_{t\leq u\leq1-\lambda/n} \frac{ n^{\nu}(1-y) \left\vert Z_{n}%
-G_{n}\right\vert } {\left\{  u(1-u)\right\}  ^{\rho}} = O_{\mathbb{P}}
\left(  n^{\nu-1/2}\log n \right)  = o_{\mathbb{P}}(1).
\label{right-binomial-error}%
\end{equation}

\noindent Combining
\eqref{right-error-decomposition}--\eqref{right-binomial-error} and using
$\mathbb{P}(A_{n}^{c})\longrightarrow0$, we obtain
\begin{equation}
\sup_{t\leq u\leq1-\lambda/n} \frac{ n^{\nu} \left\vert \alpha_{n}%
(u)-B_{n}(u)\right\vert } {\left\{  u(1-u)\right\}  ^{1/2-\nu}} =
O_{\mathbb{P}}(1). \label{right-ordinary-approximation}%
\end{equation}

\noindent Relations \eqref{left-ordinary-approximation} and
\eqref{right-ordinary-approximation} prove \eqref{ordinary-approximation}.

\smallskip\noindent\textbf{Step 6: Approximation of the increment process.}
\smallskip

\noindent We finally apply the near-one approximation to the left conditional
empirical process. Let
\[
\frac{\lambda}{n}\leq s<t, \qquad x:=\frac{s}{t}, \qquad z:=1-x.
\]
On $A_{n}$,
\[
N_{n}(1-z) = N_{n}x = \frac{N_{n}s}{t} \geq\frac{\lambda N_{n}}{nt} \geq
\frac{\lambda}{2} = a.
\]
Hence
\[
z \leq1-\frac{a}{N_{n}},
\]
which is precisely the domain required in \eqref{near-one-weighted-approximation}.

\noindent Combining \eqref{increment-exact-decomposition} and
\eqref{bridge-increment-decomposition}, we obtain
\begin{align}
\alpha_{n}(s;t)-B_{n}(s;t) ={}  &  -\sqrt{\frac{N_{n}}{n}}\, \Delta_{N_{n}%
}^{(L)}(1-x)\nonumber\\
&  - \left(  \sqrt{\frac{N_{n}}{n}}-\sqrt{t} \right)  L_{n}^{\ast}(1-x) +
x(Z_{n}-G_{n}). \label{increment-error-decomposition}%
\end{align}

For the first term, on $A_{n}$,
\begin{align*}
&  \sup_{\lambda/n\leq s<t} \frac{ n^{\nu}\sqrt{N_{n}/n}\, \left\vert
\Delta_{N_{n}}^{(L)}(1-s/t) \right\vert } {s^{\rho}}\\
&  \qquad\leq\left(  \frac{N_{n}}{n} \right)  ^{\rho} \sup_{0\leq
z\leq1-a/N_{n}} \frac{ N_{n}^{\nu} \left\vert \Delta_{N_{n}}^{(L)}(z)
\right\vert } {(1-z)^{\rho}} \sup_{0<x\leq1} \frac{x^{\rho}}{(tx)^{\rho}}.
\end{align*}
The last supremum equals $t^{-\rho}$. Thus, by the random-index observation
and \eqref{near-one-weighted-approximation},
\begin{equation}
\sup_{\lambda/n\leq s<t} \frac{ n^{\nu}\sqrt{N_{n}/n}\, \left\vert
\Delta_{N_{n}}^{(L)}(1-s/t) \right\vert } {s^{\rho}} = O_{\mathbb{P}}(1).
\label{increment-leading-error}%
\end{equation}

\noindent Using \eqref{left-square-root-error},
\[
\left\Vert L_{n}^{\ast}\right\Vert _{\infty} = O_{\mathbb{P}}(1), \qquad
s\geq\frac{\lambda}{n},
\]
and $\nu+\rho=1/2$, we obtain
\begin{equation}
\sup_{\lambda/n\leq s<t} \frac{ n^{\nu} \left\vert \sqrt{N_{n}/n}-\sqrt{t}
\right\vert \left\vert L_{n}^{\ast}(1-s/t) \right\vert } {s^{\rho}} =
O_{\mathbb{P}}(1). \label{increment-coefficient-error}%
\end{equation}

\noindent Finally,
\[
\frac{x}{s^{\rho}} = \frac{s/t}{s^{\rho}} = \frac{1}{t}s^{1-\rho} \leq C_{t}.
\]
Hence, by \eqref{binomial-normal-approximation},
\begin{equation}
\sup_{\lambda/n\leq s<t} \frac{ n^{\nu}x \left\vert Z_{n}-G_{n}\right\vert }
{s^{\rho}} = O_{\mathbb{P}} \left(  n^{\nu-1/2}\log n \right)  =
o_{\mathbb{P}}(1). \label{increment-binomial-error}%
\end{equation}

\noindent Combining
\eqref{increment-error-decomposition}--\eqref{increment-binomial-error} and
using $\mathbb{P}(A_{n}^{c})\longrightarrow0$, we obtain \eqref{increment-approximation}.

\noindent This completes the proof of the ordinary empirical-process
approximation and the empirical-increment approximation. The argument was
carried out for the originally given i.i.d.\ uniform sample. The extensions of
the probability space introduce only the Gaussian coupling variables and leave
the original sample unchanged. For the fixed terminal point $t$, the same
sequence $\{B_{n}\}_{n\geq1}$ satisfies both empirical approximations, and its
construction is independent of the particular values of $\lambda$ and $\nu$.

\subsection{Uniform quantile process and its increments}

\mbox{}

\smallskip\noindent\textbf{Step 1: Approximation of the ordinary quantile
process.} \smallskip

\noindent We first record the ordinary quantile-process approximation in the
form in which it is needed here. The weighted construction of \cite{CsCsHM86}
is a joint empirical--quantile construction: the uniform empirical process and
the uniform quantile process are approximated by the same sequence of Brownian
bridges, up to the usual sign convention for the quantile process. Since we
use
\[
\beta_{n}(u)=\sqrt n\{u-V_{n}(u)\},
\]
rather than $\sqrt n\{V_{n}(u)-u\}$, the sign is reversed and the
approximating bridge is precisely $B_{n}$. Hence, on the probability space
used above, the same sequence $\{B_{n}\}_{n\geq1}$ may be chosen so that, for
every fixed $\lambda>0$ and every $0\leq\eta<1/2$,
\begin{equation}
\sup_{\lambda/n\leq u\leq1-\lambda/n} \frac{ n^{\eta} \left\vert \beta
_{n}(u)-B_{n}(u) \right\vert } {\left\{  u(1-u)\right\}  ^{\kappa}} =
O_{\mathbb{P}}(1), \qquad\kappa=\frac12-\eta.
\label{ordinary-quantile-coupling}%
\end{equation}
This is exactly \eqref{ordinary-quantile-approximation}. No new Brownian
bridge is introduced for the ordinary quantile process; the only additional
bridge in the quantile part will be the bridge attached below to the reversed
quantile increment process. The coupling may be fixed independently of the
particular values of $\lambda$ and $\eta$. \smallskip

\noindent\textbf{Step 2: R\'enyi representation and reversal of the relevant
exponential spacings.} \smallskip

\noindent For each $n\geq1$, put
\[
N:=n+1 \qquad\text{and}\qquad j_{n}:=\left\lceil nt\right\rceil .
\]
Let
\[
U_{1:n}\leq\cdots\leq U_{n:n}
\]
denote the order statistics associated with $U_{1},\ldots,U_{n}$. By the
R\'enyi representation, the random vector
\[
\left(  U_{1:n},\ldots,U_{n:n} \right)
\]
has the same distribution as
\[
\left(  \frac{S_{1,n}}{S_{N,n}}, \ldots, \frac{S_{n,n}}{S_{N,n}} \right)  ,
\]
where
\[
S_{k,n}:=\sum_{r=1}^{k}E_{r,n}, \qquad1\leq k\leq N,
\]
and $E_{1,n},\ldots,E_{N,n}$ are i.i.d.\ standard exponential random variables.

\noindent As in the preceding constructions, the relevant spaces are standard
Borel spaces. Hence, by disintegration and randomization, after a further
extension of the probability space if necessary, we may attach to the original
order statistics, for every $n$, a vector
\[
E_{1,n},\ldots,E_{N,n}
\]
whose components are i.i.d.\ standard exponential random variables and such
that
\begin{equation}
U_{k:n} = \frac{S_{k,n}}{S_{N,n}}, \qquad k=1,\ldots,n, \quad\text{a.s.}
\label{renyi-original-order-statistics}%
\end{equation}
The original uniform sample is left unchanged by this extension.

\noindent Reverse the first $j_{n}$ exponential variables and leave the
remaining ones in their original order. More precisely, define
\[
\xi_{r,n} :=
\begin{cases}
E_{j_{n}-r+1,n}, & 1\leq r\leq j_{n},\\
E_{r,n}, & j_{n}<r\leq N.
\end{cases}
\]
Since this is a deterministic permutation of $E_{1,n},\ldots,E_{N,n}$, the
variables
\[
\xi_{1,n},\ldots,\xi_{N,n}
\]
are again i.i.d.\ standard exponential random variables. Put
\[
T_{k,n} := \sum_{r=1}^{k}\xi_{r,n}, \qquad0\leq k\leq N,
\]
with $T_{0,n}:=0$. In particular,
\[
T_{N,n}=S_{N,n}.
\]

\noindent Define
\[
U_{k:n}^{\mathrm{rev}} := \frac{T_{k,n}}{T_{N,n}}, \qquad k=1,\ldots,n.
\]
The vector
\[
\left(  U_{1:n}^{\mathrm{rev}}, \ldots, U_{n:n}^{\mathrm{rev}} \right)
\]
therefore has the distribution of the order statistics of an i.i.d.\ uniform
sample of size $n$. Let $V_{n,t}^{\mathrm{rev}}$ denote the corresponding
empirical quantile function, defined by
\[
V_{n,t}^{\mathrm{rev}}(0):=0
\]
and
\[
V_{n,t}^{\mathrm{rev}}(u) := U_{k:n}^{\mathrm{rev}} \quad\text{whenever}%
\quad\frac{k-1}{n}<u\leq\frac{k}{n}, \qquad k=1,\ldots,n.
\]
Its uniform quantile process is
\[
\beta_{n,t}^{\mathrm{rev}}(u) := \sqrt{n}\left\{  u-V_{n,t}^{\mathrm{rev}}(u)
\right\}  , \qquad0\leq u\leq1.
\]

\noindent For $0\leq s<t$, put
\[
\ell_{n}(s) := \left\lceil n(t-s)\right\rceil \qquad\text{and}\qquad k_{n}(s)
:= j_{n}-\ell_{n}(s).
\]
Then $0\leq k_{n}(s)\leq j_{n}-1$, and the definition of the empirical
quantile function gives
\[
V_{n}(t)=U_{j_{n}:n} \qquad\text{and}\qquad V_{n}(t-s)=U_{\ell_{n}(s):n}.
\]
By \eqref{renyi-original-order-statistics} and the reversal of the first
$j_{n}$ spacings,
\begin{align*}
V_{n}(t)-V_{n}(t-s)  &  = \frac{ E_{\ell_{n}(s)+1,n} +\cdots+ E_{j_{n},n} }
{S_{N,n}}\\
&  = \frac{T_{k_{n}(s),n}}{T_{N,n}}\\
&  = V_{n,t}^{\mathrm{rev}} \left(  \frac{k_{n}(s)}{n} \right)  ,
\qquad\text{a.s.}%
\end{align*}
where the last identity also holds when $k_{n}(s)=0$. Consequently,
\begin{equation}
\beta_{n}(s;t) = \beta_{n,t}^{\mathrm{rev}} \left(  \frac{k_{n}(s)}{n}
\right)  + \sqrt{n}\left\{  s-\frac{k_{n}(s)}{n} \right\}  .
\label{quantile-increment-reversal-identity}%
\end{equation}

\noindent Finally,
\[
k_{n}(s)-ns = \left\{  \left\lceil nt\right\rceil -nt \right\}  - \left\{
\left\lceil n(t-s)\right\rceil -n(t-s) \right\}  ,
\]
and therefore
\begin{equation}
\left\vert k_{n}(s)-ns \right\vert <1. \label{quantile-index-rounding}%
\end{equation}
In particular,
\[
\left\vert s-\frac{k_{n}(s)}{n} \right\vert < \frac{1}{n}.
\]

\newpage\noindent\textbf{Step 3: Coupling of the reversed quantile process and
completion of the proof.} \smallskip

\noindent We now attach a Brownian bridge to the reversed quantile process
constructed in Step~2. This bridge is used only as an auxiliary object in the
proof. The theorem itself is stated in terms of the increment $B_{n}%
(s;t)=B_{n}(t)-B_{n}(t-s)$ of a standard Brownian bridge; this is legitimate
because the auxiliary bridge constructed below has exactly the same law as
this increment process. By the same disintegration and randomization argument
used previously, after enlarging the probability space if necessary, we may
construct a sequence of standard Brownian bridges
\[
\left\{  B_{n,t}^{\left(  q\right)  }(u);\,0\leq u\leq1 \right\}  _{n\geq1}
\]
such that the joint distribution of
\[
\left(  \beta_{n,t}^{\mathrm{rev}}, B_{n,t}^{\left(  q\right)  } \right)
\]
is the one supplied by the weighted quantile-process coupling of
\cite{CsCsHM86}. In particular, applying that approximation with lower cutoff
$1/n$, we have
\begin{equation}
\sup_{1/n\leq u\leq1-1/n} \frac{ n^{\eta} \left\vert \beta_{n,t}%
^{\mathrm{rev}}(u) - B_{n,t}^{\left(  q\right)  }(u) \right\vert } {\left\{
u(1-u)\right\}  ^{\kappa}} = O_{\mathbb{P}}(1), \qquad\kappa=\frac{1}{2}-\eta.
\label{reversed-quantile-coupling}%
\end{equation}
The coupling may be fixed independently of the particular values of $\lambda$
and $\eta$.

\noindent We shall also use the following weighted discretization bound:
\begin{equation}
\sup_{2/n\leq s<t} \frac{ n^{\eta} \left\vert B_{n,t}^{\left(  q\right)  }
\left(  \frac{k_{n}(s)}{n} \right)  - B_{n,t}^{\left(  q\right)  }(s)
\right\vert } {s^{\kappa}} = O_{\mathbb{P}}(1).
\label{weighted-bridge-discretization}%
\end{equation}
To verify it, note first from \eqref{quantile-index-rounding} that
\[
\left\vert \frac{k_{n}(s)}{n}-s \right\vert < \frac{1}{n}.
\]
Since each $B_{n,t}^{\left(  q\right)  }$ is a standard Brownian bridge, it
has the same distribution as
\[
B(u)=W(u)-uW(1), \qquad0\leq u\leq1,
\]
where $W$ is a standard Wiener process. The contribution of the linear term is
bounded by
\[
\sup_{2/n\leq s<t} \frac{ n^{\eta-1}\left\vert W(1)\right\vert } {s^{\kappa}}
\leq2^{-\kappa}n^{-1/2}\left\vert W(1)\right\vert = o_{\mathbb{P}}(1).
\]

\noindent It remains to control the Wiener-process increments. Divide
$[2/n,t]$ into the dyadic blocks
\[
I_{j,n} := \left[  \frac{2^{j}}{n}, \frac{2^{j+1}}{n} \right]  , \qquad
j=1,\ldots,J_{n},
\]
where $J_{n}$ is the largest integer for which $2^{J_{n}}/n<t$. On $I_{j,n}$,
\[
s^{-\kappa} \leq n^{\kappa}2^{-j\kappa},
\]
and hence, since $\eta+\kappa=1/2$,
\begin{align*}
&  \sup_{s\in I_{j,n}} \frac{ n^{\eta} \left\vert W\left(  k_{n}(s)/n\right)
-W(s) \right\vert } {s^{\kappa}}\\
&  \qquad\leq\frac{\sqrt{n}}{2^{j\kappa}} \sup_{\substack{0\leq u,v\leq
t+1/n:\\\left\vert u-v\right\vert \leq1/n,\\u\vee v\in I_{j,n}^{+}}}
\left\vert W(u)-W(v)\right\vert ,
\end{align*}
where
\[
I_{j,n}^{+} := \left[  \frac{2^{j}-1}{n}, \frac{2^{j+1}+1}{n} \right]  .
\]
By Brownian scaling and the usual Gaussian maximal inequality, there exist
constants $C,c>0$, independent of $j$ and $n$, such that, for every $M>0$,
\[
\mathbb{P}\left\{  \sup_{s\in I_{j,n}} \frac{ n^{\eta} \left\vert W\left(
k_{n}(s)/n\right)  -W(s) \right\vert } {s^{\kappa}} >M \right\}  \leq C2^{j}
\exp\left\{  -cM^{2}2^{2j\kappa} \right\}  .
\]
Since $\kappa>0$,
\[
\sum_{j=1}^{\infty} 2^{j} \exp\left\{  -cM^{2}2^{2j\kappa} \right\}
\longrightarrow0 \qquad\text{as}\qquad M\longrightarrow\infty.
\]
This proves \eqref{weighted-bridge-discretization}.

\noindent We now combine the preceding bounds with the exact identity
\eqref{quantile-increment-reversal-identity}. We first consider
\[
\frac{2}{n}\leq s<t.
\]
By \eqref{quantile-index-rounding},
\[
\frac{k_{n}(s)}{n} \geq s-\frac{1}{n} \geq\frac{s}{2}
\]
and
\[
\frac{k_{n}(s)}{n} \leq s+\frac{1}{n} \leq\frac{3s}{2}.
\]
Moreover, since $t<1$, for all sufficiently large $n$,
\[
\frac{k_{n}(s)}{n} \leq t+\frac{1}{n} \leq\frac{1+t}{2} < 1.
\]
Consequently, \eqref{reversed-quantile-coupling} implies
\begin{equation}
\sup_{2/n\leq s<t} \frac{ n^{\eta} \left\vert \beta_{n,t}^{\mathrm{rev}}
\left(  k_{n}(s)/n \right)  - B_{n,t}^{\left(  q\right)  } \left(  k_{n}(s)/n
\right)  \right\vert } {s^{\kappa}} = O_{\mathbb{P}}(1).
\label{reversed-quantile-grid-error}%
\end{equation}
Indeed, the ratio
\[
\frac{ \left[  \left\{  k_{n}(s)/n\right\}  \left\{  1-k_{n}(s)/n\right\}
\right]  ^{\kappa} } {s^{\kappa}}
\]
is uniformly bounded over $2/n\leq s<t$.

\noindent The deterministic rounding term in
\eqref{quantile-increment-reversal-identity} satisfies
\begin{equation}
\sup_{\lambda/n\leq s<t} \frac{ n^{\eta}\sqrt{n} \left\vert s-\frac{k_{n}%
(s)}{n} \right\vert } {s^{\kappa}} \leq\sup_{\lambda/n\leq s<t} \frac
{n^{\eta-1/2}}{s^{\kappa}} \leq\lambda^{-\kappa}.
\label{quantile-rounding-error}%
\end{equation}
It follows from \eqref{quantile-increment-reversal-identity},
\eqref{weighted-bridge-discretization}, \eqref{reversed-quantile-grid-error},
and \eqref{quantile-rounding-error} that
\begin{equation}
\sup_{\max(\lambda,2)/n\leq s<t} \frac{ n^{\eta} \left\vert \beta
_{n}(s;t)-B_{n,t}^{\left(  q\right)  }(s) \right\vert } {s^{\kappa}} =
O_{\mathbb{P}}(1). \label{quantile-increment-away-from-zero}%
\end{equation}

\noindent It remains only to treat the finite-count interval
\[
\frac{\lambda}{n} \leq s < \frac{2}{n},
\]
which is nonempty only when $\lambda<2$. On this interval,
\eqref{quantile-index-rounding} implies
\[
0\leq k_{n}(s)\leq2.
\]
Hence, by the reversed R\'enyi representation,
\[
0 \leq V_{n}(t)-V_{n}(t-s) = \frac{T_{k_{n}(s),n}}{T_{N,n}} \leq\frac{T_{2,n}%
}{T_{N,n}}.
\]
Since
\[
T_{2,n} = O_{\mathbb{P}}(1) \qquad\text{and}\qquad\frac{T_{N,n}}{n}
\longrightarrow1 \quad\text{in probability},
\]
we obtain
\[
\sup_{\lambda/n\leq s<2/n} \left\vert \beta_{n}(s;t) \right\vert =
O_{\mathbb{P}}\left(  n^{-1/2}\right)  .
\]
Furthermore, by the Brownian representation and Brownian scaling,
\[
\sup_{0\leq s\leq2/n} \left\vert B_{n,t}^{\left(  q\right)  }(s) \right\vert =
O_{\mathbb{P}}\left(  n^{-1/2}\right)  .
\]
Since $s\geq\lambda/n$ and $\eta+\kappa=1/2$, it follows that
\begin{equation}
\sup_{\lambda/n\leq s<2/n} \frac{ n^{\eta} \left\vert \beta_{n}(s;t)-B_{n,t}%
^{\left(  q\right)  }(s) \right\vert } {s^{\kappa}} = O_{\mathbb{P}}(1).
\label{quantile-increment-finite-count-range}%
\end{equation}

\noindent Combining \eqref{quantile-increment-away-from-zero} and
\eqref{quantile-increment-finite-count-range}, we obtain the auxiliary
coupling bound
\[
\sup_{\lambda/n\leq s<t} \frac{ n^{\eta} \left\vert \beta_{n}(s;t)-B_{n,t}%
^{\left(  q\right)  }(s) \right\vert } {s^{1/2-\eta}} = O_{\mathbb{P}}(1).
\]
It remains to identify the distribution of the Gaussian process appearing in
this approximation. For $0\leq r\leq s\leq t$,
\[
\operatorname{Cov}\left\{  B_{n,t}^{\left(  q\right)  }(r), B_{n,t}^{\left(
q\right)  }(s) \right\}  = r-rs.
\]
On the other hand, if $B$ is a standard Brownian bridge, a direct covariance
calculation gives
\[
\operatorname{Cov}\left\{  B(t)-B(t-r), B(t)-B(t-s) \right\}  = r-rs.
\]
Both processes are centered Gaussian processes with continuous sample paths.
Therefore,
\[
\left\{  B_{n,t}^{\left(  q\right)  }(s);\,0\leq s\leq t \right\}
\]
has the same distribution as
\[
\left\{  B(t)-B(t-s);\,0\leq s\leq t \right\}  .
\]
Consequently, for the weak-approximation formulation used in
Theorem~\ref{Theo1}, the auxiliary bridge may be represented by the increment
$B_{n}(s;t)=B_{n}(t)-B_{n}(t-s)$ of a standard Brownian bridge. This gives
\eqref{quantile-increment-approximation} as stated.

\smallskip

\noindent The ordinary quantile-process approximation and the
quantile-increment approximation are now proved. Together with the two
empirical approximations established in the preceding subsection, this
completes the proof of Theorem~\ref{Theo1}. All extensions of the probability
space leave the original uniform sample unchanged, and the auxiliary
constructions may be chosen independently of the particular values of
$\lambda$, $\nu$, and $\eta$.

\renewcommand{\refname}{References}


\begin{thebibliography}{99}                                                                                               %


\bibitem {Alex87}Alexander, K.S.: The central limit theorem for weighted
empirical processes indexed by sets. \textit{J. Multivariate Anal.}
\textbf{22}, 313--339 (1987). \url{https://doi.org/10.1016/0047-259X(87)90093-5}

\bibitem {BMN2015}Brahimi, B., Meraghni, D., Necir, A.: Gaussian approximation
to the extreme value index estimator of a heavy-tailed distribution under
random censoring. \textit{Math. Methods Statist.} \textbf{24}, 266--279
(2015). \url{https://doi.org/10.3103/S106653071504002X}

\bibitem {CsCsHM86}Cs\"org\H{o}, M., Cs\"org\H{o}, S., Horv\'ath, L., Mason,
D.M.: Weighted empirical and quantile processes. \textit{Ann. Probab.}
\textbf{14}, 31--85 (1986). \url{https://doi.org/10.1214/aop/1176992617}

\bibitem {cdm85}Cs\"org\H{o}, S., Deheuvels, P., Mason, D.M.: Kernel estimates
of the tail index of a distribution. \textit{Ann. Statist.} \textbf{13},
1050--1077 (1985). \url{https://doi.org/10.1214/aos/1176349656}

\bibitem {CSHM1986}Cs\"org\H{o}, S., Horv\'ath, L., Mason, D.M.: What portion
of the sample makes a partial sum asymptotically stable or normal?
\textit{Probab. Theory Relat. Fields} \textbf{72}, 1--16 (1986). \url{https://doi.org/10.1007/BF00343893}

\bibitem {CSM1985}Cs\"org\H{o}, S., Mason, D.M.: Central limit theorems for
sums of extreme values. \textit{Math. Proc. Cambridge Philos. Soc.}
\textbf{98}, 547--558 (1985). \url{https://doi.org/10.1017/S0305004100063751}

\bibitem {CSM1986}Cs\"org\H{o}, S., Mason, D.M.: The asymptotic distribution
of sums of extreme values from a regularly varying distribution. \textit{Ann.
Probab.} \textbf{14}, 974--983 (1986). \url{https://doi.org/10.1214/aop/1176992451}

\bibitem {DeEn96}Deheuvels, P., Einmahl, J.H.J.: On the strong limiting
behavior of local functionals of empirical processes based upon censored data.
\textit{Ann. Probab.} \textbf{24}, 504--525 (1996). \url{https://doi.org/10.1214/aop/1042644729}

\bibitem {EnKo92}Einmahl, J.H.J., Koning, A.J.: Limit theorems for a general
weighted process under random censoring. \textit{Canad. J. Statist.}
\textbf{20}, 77--89 (1992). \url{https://doi.org/10.2307/3315576}

\bibitem {GNM2026A}Guesmia, N.E., Necir, A., Meraghni, D.: Nelson--Aalen
kernel estimator to the tail index of right censored Pareto-type data.
\textit{J. Korean Statist. Soc.} \textbf{55}, 627--667 (2026). \url{https://doi.org/10.1007/s42952-025-00362-y}

\bibitem {GNM2026B}Guesmia, N.E., Necir, A., Meraghni, D.: Adapted kernel
estimator to the tail index of randomly right-censored Pareto-type data.
\textit{Jpn. J. Stat. Data Sci.} (2026). \url{https://doi.org/10.1007/s42081-026-00336-2}

\bibitem {KMT75}Koml\'os, J., Major, P., Tusn\'ady, G.: An approximation of
partial sums of independent random variables and the sample distribution
function. I. \textit{Z. Wahrscheinlichkeitstheorie und Verw. Gebiete}
\textbf{32}, 111--131 (1975). \url{https://doi.org/10.1007/BF00533093}

\bibitem {KMT76}Koml\'os, J., Major, P., Tusn\'ady, G.: An approximation of
partial sums of independent random variables and the sample distribution
function. II. \textit{Z. Wahrscheinlichkeitstheorie und Verw. Gebiete}
\textbf{34}, 33--58 (1976). \url{https://doi.org/10.1007/BF00532688}

\bibitem {MNS2025}Meraghni, D., Necir, A., Soltane, L.: Nelson--Aalen tail
product-limit process and extreme value index estimation under random
censorship. \textit{Sankhya A} \textbf{87}, 526--574 (2025). \url{https://doi.org/10.1007/s13171-025-00384-y}

\bibitem {peng2001}Peng, L.: Estimating the mean of a heavy-tailed
distribution. \textit{Statist. Probab. Lett.} \textbf{52}, 255--264 (2001). \url{https://doi.org/10.1016/S0167-7152(00)00203-0}

\bibitem {SW1982}Shorack, G.R., Wellner, J.A.: Limit theorems and inequalities
for the uniform empirical process indexed by intervals. \textit{Ann. Probab.}
\textbf{10}, 639--652 (1982). \url{https://doi.org/10.1214/aop/1176993773}
\end{thebibliography}
\end{document}